\documentclass[11pt]{amsart}
\usepackage{amssymb,mathrsfs,graphicx,enumerate,color}
\usepackage{amsmath, amsfonts, amscd, epsfig}

\usepackage{mathtools}
\usepackage{hyperref}
\usepackage{comment}

\newcommand{\Ccal}{\mathcal{C}}
\newcommand{\Jcal}{\mathcal{J}}
\newcommand{\Ocal}{\mathcal{O}}

\renewcommand{\hbar}{\bar{h}}

\newcommand{\wbar}{\bar{w}}

\newcommand{\ptil}{\tilde{p}}

\newcommand{\util}{\tilde{u}}
\newcommand{\vtil}{\tilde{v}}
\newcommand{\wtil}{\tilde{w}}

\newcommand{\abs}[1]{\left|#1\right|}

\DeclareMathOperator{\tr}{tr}

\newcommand{\g}{\gamma}

\newcommand{\e}{\varepsilon}

\renewcommand{\th}{\theta}
\renewcommand{\k}{\kappa}
\renewcommand{\l}{\lambda}

\newcommand{\m}{\mu}

\newcommand{\x}{\xi}

\renewcommand{\r}{\rho}
\newcommand{\s}{\sigma}
\renewcommand{\t}{\tau}

\newcommand{\rd}{\partial}

\newcommand{\thtil}{\tilde{\th}}

\usepackage{colortbl}
\definecolor{black}{rgb}{0.0, 0.0, 0.0}
\definecolor{red}{rgb}{1.0, 0.5, 0.5}

\title[   ]{Traveling Wave Solutions to Brenner-Navier-Stokes-Fourier system}

\author[Eo]{Saehoon Eo}
\address[Saehoon Eo]
{ Department of Mathematical Sciences, \newline
Korea Advanced Institute of
Science and Technology \\
Daejeon 34141, Korea}
\email{eosehoon@kaist.ac.kr}

\author[Eun]{Namhyun Eun}
\address[Namhyun Eun]
{ Department of Mathematical Sciences, \newline
Korea Advanced Institute of
Science and Technology \\
Daejeon 34141, Korea}
\email{namhyuneun@kaist.ac.kr}
 
\author[Kang]{Moon-Jin Kang}
\address[Moon-Jin Kang]
{ Department of Mathematical Sciences, \newline
Korea Advanced Institute of
Science and Technology \\
Daejeon 34141, Korea}
\email{moonjinkang@kaist.ac.kr}

\author[Oh]{HyeonSeop Oh}
\address[HyeonSeop Oh]
{ Department of Mathematical Sciences, \newline
Korea Advanced Institute of
Science and Technology \\
Daejeon 34141, Korea}
\email{ohs2509@kaist.ac.kr}

\newtheorem{theorem}{Theorem}[section]
\newtheorem{lemma}{Lemma}[section]

\newtheorem{proposition}{Proposition}[section]
\newtheorem{remark}{Remark}[section]

\newtheorem{definition}{Definition}[section]

\newcommand{\beq}{\begin{equation}}
\newcommand{\eeq}{\end{equation}}
\newcommand{\bsp}{\begin{split}}
\newcommand{\esp}{\end{split}}

\newcommand{\bbr}{\mathbb R}

\newcommand{\eps}{\varepsilon}

\numberwithin{figure}{section}
\newcommand{\RR}{{\mathbb{R}}}

\newcommand{\vt}{{\tilde{v}_\e}}
\newcommand{\ut}{{\tilde{u}_\e}}
\newcommand{\pt}{{\tilde{p}_\e}}
\newcommand{\tht}{{\tilde{\th}_\e}}

\newcommand{\step}[1]{\vskip0.2cm \noindent{\it Step #1:} }

\newcommand{\sbare}{\overline{\sigma}_\varepsilon}

\begin{document}

\date{\today}


\thanks{\textbf{Acknowledgment.} The authors thank Professor Alexis Vasseur for valuable comments on BNSF system.   
This work was supported by Samsung Science and Technology Foundation under Project Number SSTF-BA2102-01. }

\begin{abstract}
As a continuum model for compressible fluid flows, Howard Brenner proposed the so-called Brenner-Navier-Stokes-Fourier(BNSF) system that improves some flaws of the Navier-Stokes-Fourier(NSF) system. For BNSF system, the volume velocity concept is introduced and is far different from the mass velocity of NSF, since the density of a compressible fluid is inhomogeneous. Although BNSF was introduced more than ten years ago, the mathematical study on BNSF is still in its infancy. 
We consider the BNSF system in the Lagrangian mass coordinates.
We prove the existence and uniqueness of monotone traveling wave solutions to the BNSF system. We also present some quantitative estimates for them.   
\end{abstract}
\maketitle \centerline{\date}

\tableofcontents

\section{Introduction}
\setcounter{equation}{0}
We consider the so-called Brenner-Navier-Stokes-Fourier(BNSF) system (in the Eulerian coordinates) in 1D :
\begin{align} \label{BNSF-E}
\left\{
\begin{aligned}
&\rd_t \r + \rd_x (\r u_m) = 0, \\
&\rd_t (\r u_v)+\rd_x (\r u_v u_m + p)=\rd_x \left(\m\rd_x u_v\right), \\ 
&\rd_t \left(\r \left(\frac{u_v^2}{2} + e \right) \right) + \rd_x \left(\r\left( \frac{u_v^2}{2} + e \right)u_m+pu_v\right)
= \rd_x \left(\m u_v \rd_x u_v + \k \rd_x \th\right).
\end{aligned}
\right.
\end{align}
 Here, \(\r\) is the density, \(u_m\) is the mass velocity, \(u_v\) is the volume velocity, \(p\) is the pressure, and \(e\) is the internal energy. In addition, \(\m\) and \(\k\) represent the viscosity coefficient and the heat-conductivity coefficient, respectively, which are some positive constants depending only on the type of gas.

In a series of papers \cite{brenner2005kinematics, brenner2005navier, brenner2006fluid, brenner2012beyond}, Howard Brenner proposed the system \eqref{BNSF-E} to improve some defects of the Navier-Stokes-Fourier(NSF) system. 
His idea is based on so-called `bi-velocity theory' which indicates the existence of two different velocities: one is the mass velocity \(u_m\) which is the classical concept; the other is the volume velocity \(u_v\). 
Note that these two velocities are implicitly assumed to be identical in the vast majority of studies on continuum fluid mechanics. This perspective has remained unchallenged since Euler's era. 
However, Brenner contended that, in general, \(u_m \neq u_v\), and the inconsistency increases as the density gradient becomes larger. In fact, there have been numerous studies pointing out and attempting to improve the imperfections of the NSF system(see \cite{arkilic2001mass, chakraborty2007derivations, dadzie2010volume, dadzie2012analysis, dadzie2012spatial, dadzie2008continuum, dongari2010predicting, dongari2009extended, durst2011treatments, gad1999fluid, goddard2010material, greenshields2007structure, harley1995gas, klimontovich1992need, klimontovich1993hamiltonian, klimontovich1995unified, koide2012navier, Koide2011NavierStokesEB, reese2003new}), but it can be said that Brenner's approach is the most systematic and well-established.

The necessity of the concept of a new velocity can be explained as follows: consider the work generated by the displacement of gas particles due to the pressure within the gas. To calculate the amount of the work, the displacement distance must be defined. However, it is crucial to point out that it is not the individual gas particle that forms pressure and is displaced by pressure, but rather a collection of gas particles. Since the boundary of a collection of gas particles cannot be clearly defined, it intuitively seems impossible to rigorously define the displacement distance. This is related to the fact that the volume of gas is not a mass-point property. The volume of gas cannot be understood as a property possessed by particles at a point but can be considered statistically as a property of a collection  of particles at best. Hence,  the volume velocity, as a new velocity, is required to describe the motion of a collection of particles, as opposed to the conventional velocity. Therefore, in the existing NSF system, if we only use the mass velocity for the mass transportation and the convection terms, but replace all other velocities by the volume velocity, it results in the BNSF system \eqref{BNSF-E}. Note that the concept of volume velocity  is inherently elusive. Thus, it must be defined by a certain constitutive equation.

Brenner proposed the volume velocity concept through the following argument. 
First, he derived the constitutive equations from linear irreversible thermodynamics. 
He then demonstrated that these equations align with the Burnett's solution to the Boltzmann equation for a dilute monatomic case. 
Additionally, it turns out that they better match the experimental data. 
One remark is that when the no-slip boundary condition at solid surfaces is described by the volume velocity, it more closely corresponds to the experimental data as well. 
(See \cite{brenner2012beyond})

The constitutive relation between $u_v$ and $u_m$ suggested by Brenner  is as follows:
\begin{equation} \label{umuv}
u_v=u_m + \frac{\k}{\r c_p} \nabla \ln \r,
\end{equation}
where \(c_p\) is the specific heat at constant pressure. The difference between the two velocities is proportional to \(\nabla \r\). Thus, in incompressible fluids with uniform density, the concept of volume velocity becomes insignificant. Brenner mentioned that the NSF system is a PDE limited to incompressible fluids, and it seems it took a long time to discover the flaws of the NSF system because, during its initial validation, data from liquids or essentially incompressible gases were mainly used. (See \cite{brenner2012beyond})

Although this amendment to the NSF system were proposed over ten years ago, the mathematical study on BNSF is still in its infancy. We refer to Feireisl-Vassuer \cite{FeVa}, which proved the existence of the global weak solutions to the initial boundary value problem.\\

In this paper, we consider the one-dimensional BNSF system \eqref{BNSF-E} in the Lagrangian mass coordinates:
(still use the same notation $x$ to denote the mass variable $y$)
\begin{equation}
\left\{
\begin{aligned} \label{main}
    &v_t-u_x = \left(\t\frac{v_x}{v}\right)_x, \\
    &u_t+p(v,\th)_x = \left(\m\frac{u_x}{v}\right)_x, \\
    &E_t+(p(v,\th)u)_x = 
    \left(\k\frac{\th_x}{v}\right)_x+\left(\m\frac{uu_x}{v}\right)_x,
\end{aligned} \right.
\end{equation}
where \(v\) is the specific volume, \(u\) is now the volume velocity, \(\th\) is the absolute temperature, \(E=e+\frac{u^2}{2}\) is the total energy, and \(\t = \k/c_p\) is the Brenner coefficient. When converting to Lagrangian mass coordinates with the constitutive relation \eqref{umuv}, the convection terms disappear, thus Brenner's modification only remains in the mass conservation law. We also consider the ideal polytropic gas where the pressure law \(p\) and the internal energy \(e\) function are given by    
\begin{equation} \label{pressure}
p(v,\th)=\frac{R\th}{v}, \qquad e(\th)=\frac{R}{\g-1}\th + const.
\end{equation}
where \(R>0\) is the gas constant and \(\g>1\) is the adiabatic constant.

In this article, we prove that the system \eqref{main} admits viscous shock waves connecting two end states \((v_-,u_-,E_-)\) and \((v_+,u_+,E_+)\) when these two end states are close enough and satisfy the Rankine-Hugoniot condition and Lax entropy condition as follows:
\begin{equation}
\begin{aligned} \label{end-con}
&\exists~\s \quad \text{s.t.}~\left\{
\begin{aligned}
&-\s(v_+-v_-) -(u_+-u_-) =0, \\
&-\s(u_+-u_-) +p(v_+,\th_+)-p(v_-,\th_-) =0, \\
&-\s(E_+-E_-) +p(v_+,\th_+)u_+ -p(v_-,\th_-)u_- =0,
\end{aligned} \right. \\
&\text{and either \(v_->v_+\), \(u_->u_+\) and \(\th_-<\th_+\) or \(v_-<v_+\), \(u_->u_+\) and \(\th_->\th_+\) holds.}
\end{aligned}
\end{equation}
Here, \(\th_-\) and \(\th_+\) satisfy 
\[
E_- = \frac{R}{\g-1}\th_- + \frac{u_-^2}{2}, \qquad
E_+ = \frac{R}{\g-1}\th_+ + \frac{u_+^2}{2}.
\]
In other words, for given sufficiently close constant states \((v_-,u_-,E_-)\) and \((v_+,u_+,E_+)\) satisfying \eqref{end-con}, we prove the existence of a viscous shock wave (\(\vtil,\util,\thtil)(\x)=(\vtil, \util, \thtil\))(\(x-\s t\)), as a traveling wave solution to the following system of ODEs:
\begin{equation}
\left\{
\begin{aligned} \label{shock_0}
    &-\s \vtil'-\util' = \left(\t\frac{\vtil'}{\vtil}\right)', \\
    &-\s \util'+p(\vtil,\thtil)' = \left(\mu\frac{\util'}{\vtil}\right)', \\
    &-\s \left(\frac{R}{\g-1}\thtil + \frac{\util^2}{2}\right)'+(p(\vtil,\thtil)\util)' = 
    \Bigg(\k\frac{\thtil'}{\vtil}\Bigg)'+ \Bigg(\m\frac{\util \util'}{\vtil}\Bigg)', \\
    &\lim_{\x\to\pm\infty} (\vtil,\util,\thtil)(\x)=(v_\pm,u_\pm,\th_\pm).
\end{aligned}
\right. \qquad\left('=\frac{d}{d\x}\right)
\end{equation}
Here, if \(v_->v_+\), (\(\vtil, \util, \thtil\))(\(x-\s t\)) is a 1-shock wave with velocity \(\s=- \sqrt{-\frac{p_+-p_-}{v_+-v_-}}\), where \(p_\pm := p(v_\pm, \th_\pm)\).
If \(v_-<v_+\), it becomes a 3-shock wave with \(\s= \sqrt{-\frac{p_+-p_-}{v_+-v_-}}\).

\subsection{Main results}
Our first result is for the existence and uniqueness of traveling wave solutions of \eqref{main} as monotone profiles satisfying \eqref{shock_0}.

\begin{theorem}\label{thm_EU} (Existence and Uniqueness) For a given left-end state \((v_-,u_-,\th_-) \in \RR^+ \times \RR \times \RR^+\), there exists a constant \(\e_0>0\) such that for any right-end state \((v_+, u_+, \th_+) \in \RR^+ \times \RR \times \RR^+\) satisfying \eqref{end-con} and \(\e := \abs{v_--v_+} < \e_0\), there is a unique traveling wave solution \((\vt, \ut, \tht)\colon \RR\to \RR^+ \times \RR \times \RR^+\) to \eqref{main}, as a monotone profile satisfying \eqref{shock_0}, which more precisely satisfies \eqref{1shock} when $v_->v_+$, or \eqref{derivsign} when $v_-<v_+$.
 
\end{theorem}

The second result provides quantitative estimates of traveling wave solutions. 
Especially, we present estimates for the ratio between the components $\vtil, \util, \thtil$ of the traveling wave, and estimates for the exponential tail of the wave. This could be useful in stability estimates for traveling waves. 

\begin{theorem}\label{thm_estimates} (Estimates)
For a given left-end state \((v_-,u_-,\th_-) \in \RR^+ \times \RR \times \RR^+\), there exist positive constants \(\e_0\), \(C\), and \(C_1\) such that for any right-end state \((v_+, u_+, \th_+) \in \RR^+ \times \RR \times \RR^+\) satisfying \eqref{end-con} and \(\e := \abs{v_--v_+} < \e_0\), the following holds. \\
Let \((\vt, \ut, \tht)\) be the monotone solution to \eqref{shock_0} with \(\vt(0) = \frac{v_-+v_+}{2}\).\\
 Then,
\begin{align}
&\begin{aligned}\label{decay}
&\abs{(\vt(\x)-v_-, \ut(\x)-u_-, \tht(\x)-\th_-)} \le C\e e^{-C_1\e\abs{\x}} \text{ for all } \x \le 0, \\
&\abs{(\vt(\x)-v_+, \ut(\x)-u_+, \tht(\x)-\th_+)} \le C\e e^{-C_1\e\abs{\x}} \text{ for all } \x \ge 0, \\
\end{aligned} \\
&\abs{(\vtil'_\e(\x), \util'_\e(\x), \thtil'_\e(\x))} \le C\e^2e^{-C_1\e\abs{\x}} \text{ for all } \x\in \RR, \label{derivdecay}\\
&\abs{(\vtil''_\e(\x), \util''_\e(\x), \thtil''_\e(\x))} \le C\e\abs{(\vtil'_\e(\x), \util'_\e(\x), \thtil'_\e(\x))} \text{ for all } \x\in \RR.
\label{2bound1}
\end{align}
It also holds that \(\abs{\vt'} \sim \abs{\ut'} \sim \abs{\tht'}\) for all \(\x\in\RR\). More explicitly, we have the following:
\begin{align}
\abs{\util'_\e(\x) + \s_*\vtil'_\e(\x)} &\le C\e\abs{\vtil'_\e(\x)} 
\label{ratio-vu} \\
\abs{\thtil'_\e(\x) + \frac{(\g-1)p_-}{R}\vtil'_\e(\x)} &\le C\e\abs{\vtil'_\e(\x)}
\label{ratio-vth}
\end{align}
which satisfies
\begin{equation} \label{sm1}
\abs{\s_\e-\s_*} \le C \e,
\end{equation}
where
\[
\s = \s_\e := -\sqrt{-\frac{p_+-p_-}{v_+-v_-}} \, \,  \text{and } \s_* := -\sqrt{\frac{\g p_-}{v_-}} = -\frac{\sqrt{\g R \th_-}}{v_-^2},
\]
or
\[
\s = \s_\e := \sqrt{-\frac{p_+-p_-}{v_+-v_-}} \, \,  \text{and } \s_* := \sqrt{\frac{\g p_-}{v_-}} = \frac{\sqrt{\g R \th_-}}{v_-^2}.
\]
In addition, if we consider the 1-shock, i.e., \(\s_\e<0\), then we have
\beq\label{1shock}
\vtil'_\e <0, \quad \util'_\e <0, \quad \thtil'_\e >0 \text{ for all } \x \in \RR,
\eeq
and for the 3-shock case, i.e., \(\s_\e>0\),
\begin{equation}\label{derivsign}
\vtil'_\e >0, \quad \util'_\e <0, \quad \thtil'_\e <0 \text{ for all } \x \in \RR.
\end{equation}
Furthermore, the following estimate holds:
\begin{equation} \label{tail-0}
C^{-1}(v_+-\vt)(\vt-v_-) \le \vt' \le C(v_+-\vt)(\vt-v_-).
\end{equation}
\end{theorem}

Note that it suffices to prove the Theorem \ref{thm_estimates} for 3-shocks. Indeed, the results for the 1-shocks can be obatained by the change of variables \(x \mapsto -x, u \mapsto -u, \s_\e \mapsto -\s_\e\) and \(\s_* \mapsto -\s_*\). Thus, in the sequel, we only consider the 3-shock case, i.e.,
\[
\s_\e = \sqrt{-\frac{p_+-p_-}{v_+-v_-}}>0 \, \, \text{and }\s_* \coloneqq \sqrt{\frac{\g p_-}{v_-}} = \frac{\sqrt{\g R \th_-}}{v_-^2}
\]
\begin{remark}
The estimates of Theorem \ref{thm_estimates} would be useful in the stability estimates of traveling wave solutions, based on the method of $a$-contraction with shifts, as in previous results \cite{CKKV,HKKL,KV21,KV-Inven,KV-2shock,KVW23,KVW-NSF}. 
Recently, the contraction (up to shift) of large perturbations from the traveling waves was proved in \cite{EEK}, for which the estimates of Theorem \ref{thm_estimates} is crucially used.
\end{remark}

\section{Main ideas and methodology} \label{sec:pre}
\setcounter{equation}{0}

\subsection{Main ideas for the proofs}
First of all, we may integrate the ODE system \eqref{shock_0} over $(-\infty,\x]$ and use the second equation of it, to find the system of ODEs of first order:
\beq\label{temsys}
    \left\{
    \begin{aligned}
        & -\s_\e(\vtil-v_-) - (\util-u_-) = \t \frac{\vtil'}{\vtil}, \\
        & -\s_\e(\util-u_-) + (\ptil-p_-) = \m \frac{\util'}{\vtil}, \\
        &  -\s_\e\left(\frac{R}{\g-1}(\thtil-\th_-) + \frac{1}{2}(\util^2-u_-^2)\right)
    + (\ptil\util-p_-u_-) = \k \frac{\thtil'}{\vtil}.
    \end{aligned}
    \right.
\eeq
For the existence theory of solutions to the system \eqref{temsys}, it might be enough to use some typical approach such as Cauchy-Lipschitz theory. 
However, in order to get the ratio estimates as in \eqref{ratio-vu}-\eqref{ratio-vth}, we need to first show that the first derivative of wave dominates its second derivative as in \eqref{2bound1}, since we may verify \eqref{ratio-vu} from the second order equation $\eqref{shock_0}_1$, and  \eqref{ratio-vth} from the following system of two equations on $\vtil, \thtil$ variables only:
\begin{equation} \label{vthODE}
\left\{
\begin{aligned}
    & \s_\e^2(\vtil-v_-) +\t\s_\e\frac{\vtil'}{\vtil} + (\ptil - p_-)
    = \m \left(-\s_\e \frac{\vtil'}{\vtil} - \t \frac{\vtil''}{\vtil^2} + \t \frac{(\vtil')^2}{\vtil^3} \right), \\
    &-\s_\e\left(\frac{R}{\g-1}(\thtil-\th_-) - \frac{1}{2}\left(\s_\e(\vtil-v_-)+ \t \frac{\vtil'}{\vtil}\right)^2\right)
    - p_-\left(\s_\e(\vtil-v_-) + \t \frac{\vtil'}{\vtil}\right) = \k \frac{\thtil'}{\vtil},
\end{aligned}
\right.
\end{equation}
where the above system can be derived from \eqref{temsys}, and the second order term $\vtil''$ is due to the right-hand side of $\eqref{temsys}_1$. This is very different from the case for the Navier-Stokes-Fourier system that has a simpler structure as in \cite{KVW-NSF}. \\

Since we do not know the monotonicity of solutions prior to the existence of them, the proof for the desired estimates of Theorem \ref{thm_estimates} is not obvious. Thus, for the proofs of the main results, especially to verify the quantitative estimates as desired, we will employ a geometric approach to the singular perturbation system of ODEs associated to \eqref{temsys}.  To derive the singular perturbation system, we would consider the third order ODE for $\vtil$ variable only, like \eqref{veqn} that could be derived from \eqref{vthODE}.
We may (formally) introduce a new variable $v(z)$ defined by $\vtil(\x) = \e v(\e \x)+v_-$, where the scale $z:=\e \x$ is slow while $\x$ is fast.   
Then, by setting \(w_0 \coloneqq v\), \(w_1 \coloneqq v'\), \(w_2 \coloneqq \e v''\), and so \(w'_0 = w_1\), \(\e w_1' = w_2\) and \(\e w_2' = \e^2 v'''\), we may have the singular perturbation system \eqref{sys}, which is a starting point in our analysis. We may consider a locally invariant manifold near a critical manifold of \eqref{sys}, where the locally invariant manifold could be explicitly described based on Fenichel's first theorem.

\subsection{Fenichel's First Theorem and locally invariant manifolds}
We here present the main concepts for the geometric approach.
Consider the following system of ODEs with a small parameter \(\e\): 
\begin{equation}\label{slow}
\left\{
\begin{aligned}
    \e x' &= f(x, y, \e) \\
    y' &= g(x, y, \e), 
\end{aligned}
\right.
\end{equation}
where \(x\in \RR^n\), \(y\in \RR^l\). Here, we assume \(f\) and \(g\) to be smooth on a set \(U\times I\), where \(U\subset \RR^{n+l}\) is open and \(I\) is an open interval containing \(0\). 
Following \cite{JONES}, we say an \(l\)-dimensional manifold \(M_0\subset \RR^{n+l}\) is a \textit{critical manifold} of the system \eqref{slow} if each \((x, y)\in M_0\) satisfies that \(f(x, y, 0) = 0\).
To state Fenichel's first theorem, we need the following definitions:
\begin{definition}\label{def-normally_hyperbolic}
    We say that the manifold \(M_0\) is normally hyperbolic relative to \eqref{slow} if the \(n\times n\) matrix 
\[
\left. D_x f(x, y, \e) \right|_{\e=0}
\]
has \(n\) eigenvalues (counting multiplicity) with nonzero real part for each \((x, y)\in M_0\). 
\end{definition}

Fenichel's first theorem proves the existence of a manifold on which the flow of the solution of \eqref{slow} is confined. To be precise, we define the concept of a locally invariant manifold: 
\begin{definition}\label{def-locally_invariant}
    A set \(M \subseteq \RR^{n+l}\) is locally invariant under the flow from \eqref{slow} if there is an open set \(V \subseteq \RR^{n+l}\) containing \(M\) such that for any \(x_0\in M\) and any \(T\in\RR\), the trajectory \(x\) starting at \(x_0\) with \(x([0, T])\subset V\) also satisfies \(x([0, T])\subset M\). (When \(T<0\), replace \([0,T]\) with \([T,0]\)). 
\end{definition}
Note that, if \(M_0\) is normally hyperbolic, we can express it locally as a graph of \(x\) in terms of \(y\) thanks to the implicit function theorem. 
However, in the following discussion, we only deal with the case when \(M_0\) can be globally expressed as a graph of \(x\) in terms of \(y\) on a compact domain \(K\subset \RR^l\), and in this case, Fenichel's first theorem  can be stated as follows. 
\begin{proposition}  \cite{JONES} \label{prop_FFT}
    Consider the system \eqref{slow} with normally hyperbolic critical manifold \(M_0\).
    Suppose that there is a smooth function \(h\colon K\to \RR^n\) of which graph is contained in \(M_0\), i.e.,
    \[
    \{(x, y)\mid x = h(y), y\in K\} \subset M_0
    \]
    where \(K\subset \RR^l\) is a compact and simply connected smooth set. 
    Then, there exists \(\e_0>0\) such that for each \(\e\) with \(\abs{\e}<\e_0\), there is a function \(h^\e\colon K\to \RR^n\) which satisfies the followings: for each \(\e \neq 0\), the graph
    \[
    \{(x, y)\mid x = h^\e(y), y\in K\}
    \]
    is contained in an \(l\)-dimensional manifold \(M_\e\) which is locally invariant under the flow of \eqref{slow}. Moreover, \(h^\e\) can be taken to be \(C^r\) for each \(r < \infty\), jointly in \(y\) and \(\e\) and \(h=h^0\) on \(K\) where $h^0$ denotes \(h^{\eps}\) at $\eps=0$. 
\end{proposition}
\begin{remark} \label{rmk-FFT}
    Since the family of functions \(h^\e\) is regular in \(\e\) as well, a consequence of Proposition \ref{prop_FFT} is that there exists a function \(s\) with two variables \(\e\) and \(y\) which satisfies 
    \[
    h^\e(y) = h(y) + \e s(y, \e).
    \]
    Note that the function \(s\) has the same regularity with \(h^\e\); i.e. \(s\) is a \(C^r\)-function for each \(r<\infty\), jointly in \(y\) and \(\e\).
\end{remark}

\section{Proof for existence and uniqueness}
\setcounter{equation}{0}
In this section, we present the proof of Theorem \ref{thm_EU}.
Instead of proving the existence of the solution to \eqref{shock_0} directly, we consider an alternative ODE system \eqref{sys} to which Proposition \ref{prop_FFT} can be applied.
We will first show the existence and uniqueness of the solution to the ODE system \eqref{sys}, and then to the system \eqref{shock_0} through a proper transformation.

For any fixed parameter \(\e\in\bbr \) small enough (to be determined later), we consider the system of ODEs :
\begin{equation}\label{sys}
\left\{
\begin{aligned}
    w_0'(z) &= w_1(z), \\
    \e w_1'(z) &= w_2(z), \\
    \e w_2'(z) &= f(w_0(z), w_1(z), w_2(z), \e),
\end{aligned}
\right.
\end{equation}
where \(z \in \RR\) and
\begin{equation}\label{ew2'}
\begin{aligned}
    &f(w_0, w_1, w_2, \e) \\
    &\quad \coloneqq \frac{1}{\e}\left(\frac{1}{\m\t\k}\right)\underbrace{\left[-\sbare^2w_0\frac{R\sbare(\e w_0+v_-)^3}{\g-1} + R(\e w_0+v_-)^2\left(\frac{\sbare p_-}{\g-1}w_0 + \sbare p_- w_0\right)\right]}_{=:\Jcal(\e)} \\
    &\quad +\frac{R\sbare(\e w_0+v_-)^3}{(\g-1)\k}\left(-\sbare\left(\frac{\m+\t}{\m\t}\right)\frac{w_1}{\e w_0+v_-}-\frac{w_2}{(\e w_0+v_-)^2}+\e^2\frac{w_1^2}{(\e w_0+v_-)^3}\right) \\
    &\quad+ \frac{R(\e w_0+v_-)^2}{\m\t\k}\left(\hspace{-2.5pt}-\frac{\sbare^3}{2}w_0^2 - \e\sbare^2\t \frac{w_0w_1}{\e w_0+v_-} - \e^2\t^2\frac{\sbare}{2}\frac{w_1^2}{(\e w_0+v_-)^2} + p_-\t\frac{w_1}{\e w_0+v_-}\hspace{-2.5pt}\right) \\
    &\quad+ (\e w_0+v_-)\Bigg(\frac{1}{\m\t}\left((p_-- \sbare^2v_--2\e\sbare^2w_0)w_1-\sbare(\m+\t) w_2\right) \\
    &\quad \hspace{28mm}\qquad+\e^2\frac{3w_1w_2}{(\e w_0+v_-)^2}-\e^4\frac{2w_1^3}{(\e w_0+v_-)^3}\Bigg).
\end{aligned}
\end{equation}
Here, \(\sbare\) is defined by
\begin{equation}\label{s(v_+)}
    \sbare \coloneqq \frac{\sqrt{\g p_-}}{\sqrt{v_- + \frac{\g + 1}{2}\e}}.
\end{equation}
It is worth mentioning that the function \(f(w_0, w_1, w_2, \e)\) is smooth for any \(\e\) in some open interval containing \(0\), and hence so is the system \eqref{sys}. Indeed, this can be verified by the argument below. We observe that \(\Jcal(\e)\) can be rewritten as follows: 
\begin{align*}
\Jcal(\e) 
&= \frac{R\sbare}{\g-1}(\e w_0+v_-)^2w_0(\g p_--\sbare^2(\e w_0+v_-)) \\
&= \frac{R\sbare}{\g-1}v_-^2w_0(\g p_--\sbare^2v_-) 
+ \e\frac{R\sbare}{\g-1}(\e w_0^2+2w_0v_-)w_0(\g p_--\sbare^2v_-) \\
&\qquad
-\e\frac{R\sbare}{\g-1}(\e w_0+v_-)^2\sbare^2w_0^2.
\end{align*}
From \eqref{s(v_+)}, we have 
\begin{equation}\label{s(v_+)2}
\g p_--\sbare^2v_- = \g p_- \frac{\frac{\g+1}{2}\e}{v_- + \frac{\g+1}{2}\e}.
\end{equation}
Hence, it follows that 
\begin{align*}
\frac{1}{\e}\Jcal(\e)
&= \frac{R\sbare}{\g-1}v_-^2w_0\g p_-\frac{\frac{\g+1}{2}}{v_-+\frac{\g+1}{2}\e}
+ \frac{R\sbare}{\g-1}(\e w_0^2+2w_0v_-)w_0(\g p_--\sbare^2v_-) \\
&\qquad
-\frac{R\sbare}{\g-1}(\e w_0+v_-)^2\sbare^2w_0^2, 
\end{align*} 
which is a smooth function near \(\e = 0\). 

We then find critical points of the system \eqref{sys} for \(\e \neq 0\).
If \((w_0, w_1, w_2)\) is a critical point of the system, then from \eqref{ew2'} together with the relation \eqref{s(v_+)}, it holds that
\begin{align*}
0 &= \frac{\sbare(\e w_0+v_-)}{\g-1}(-\sbare^2w_0) + \frac{\sbare p_-}{\g-1}w_0
+ \sbare p_- w_0 -\e\frac{\sbare^3}{2}w_0^2 \\
&= \sbare w_0 \left[-\sbare^2\left(\frac{\e w_0 + v_-}{\g-1}+\frac{\e w_0}{2}\right) + \frac{\g p_-}{\g-1}\right] \\
&= \sbare^3 w_0 \left[-\left(\frac{\e w_0 + v_-}{\g-1}+\frac{\e w_0}{2}\right) + \frac{v_- + \e}{\g-1} + \frac{\e}{2}\right] 
= \sbare^3\e\left(\frac{1}{\g-1}+\frac{1}{2}\right)w_0(1-w_0),
\end{align*}
which proves that there are only two critical points \((0, 0 ,0), (1, 0, 0)\) for any \(\e \neq 0\).

\subsection{Proof of Theorem \ref{thm_EU}: Existence} \label{section E}
In this subsection, we prove the existence part of Theorem \ref{thm_EU}. The proof consists of two steps: we first prove that the system \eqref{sys} has a solution connecting the two critical points \((0, 0, 0)\) and \((1, 0, 0)\) by using Proposition \ref{prop_FFT}, and then we demonstrate that this gives a desired solution to the original system \eqref{shock_0}.

First of all, we calculate the critical manifold of the system \eqref{sys}; 
we have \(w_2 = 0\) and \(f(w_0, w_1, w_2, 0) = 0\). 
Equivalently, we have \(w_2 = 0\) and 
\begin{align*}
0 &= f(w_0, w_1, 0, 0) \\
&= \frac{1}{\m\t\k}\left(\frac{R\s_*}{\g-1}v_-\g p_-\frac{\g+1}{2}w_0 - \frac{R\s_*}{\g-1}v_-\g p_-w_0^2\right) \\
&\qquad
-\frac{R\g p_-v_-}{(\g-1)\k}\frac{\m+\t}{\m\t}w_1 - \frac{R\s_*v_-\g p_-}{\m\t\k}\frac{1}{2}w_0^2
+ \frac{Rv_-}{\m\t\k}p_-\t w_1
-\frac{1}{\m\t}(\g-1)v_-p_-w_1
\end{align*}
so that 
\begin{equation} \label{c-mfd}
w_1 = \frac{R\g\s_*}{R\t+R\g\m+(\g-1)^2\k}\frac{\g+1}{2}(w_0-w_0^2).
\end{equation}
Thus, the critical manifold \(M_0\) of \eqref{sys} can be parametrized by \(w_0\) as follows: 
\[
M_0 = \{(w_0, w_1, w_2)\mid w_1 = A(w_0-w_0^2), \, w_2 = 0\}
\]
which forms a curve in \(3\)-dimensional space. Here,
\[
A \coloneqq \frac{R\g\s_*}{R\t+R\g\m+(\g-1)^2\k}\frac{\g+1}{2}
\]
is a constant depending only on the type of gas.

To use Proposition \ref{prop_FFT}, it is required to verify that \(M_0\) is normally hyperbolic relative to the system \eqref{sys}.
Considering the definition of normally hyperbolicity, we need to compute
\begin{multline*}
\left. \begin{pmatrix}
    \frac{\rd w_2}{\rd w_1} & \frac{\rd w_2}{\rd w_2} \\
    \frac{\rd }{\rd w_1}f(w_0, w_1, w_2,\e) & \frac{\rd }{\rd w_2}f(w_0, w_1, w_2,\e)
\end{pmatrix}\right|_{\e = 0} \\
= \begin{pmatrix}
    0 & 1 \\
    -\frac{1}{\m\t\k}\frac{v_-p_-}{\g-1}(R\t+R\g\m+(\g-1)^2\k) &
    -\frac{\s_*v_-}{\m\t\k}\left(\frac{R}{\g-1}\m\t+(\m + \t)\k\right)
\end{pmatrix},
\end{multline*}
which is followed by \eqref{ew2'}.
Then, from simple calculations, we deduce that the matrix above has two real eigenvalues and both of them are negative, which proves that \(M_0\) is normally hyperbolic.

We now consider a closed interval \( K\) containing \([0,1]\). From Proposition \ref{prop_FFT} and Remark \ref{rmk-FFT}, for each small \(\e\) we obtain an 1-dimensional locally invariant manifold \(M_\e\) and smooth functions \(s_1\) and \(s_2\) defined on \(K\), which satisfies 
\begin{equation} \label{Me}
\{(w_0, w_1, w_2)\mid w_0 \in K, \, w_1 = A(w_0-w_0^2) + \e s_1(w_0, \e), \, w_2 = \e s_2(w_0, \e)\} \subset M_\e
\end{equation}
where the functions \(s_1\) and \(s_2\) are sufficiently regular jointly in \(w_0\) and \(\e\), in particular, \(C^3\).
\step{1} We will show that the system \eqref{sys} has a solution connecting the two critical points \((0, 0, 0)\) and \((1, 0, 0)\).
Here, we will prove the existence of a solution defined on $\bbr$ to the system \eqref{sys} passing through a point
\[
    (w_0, w_1, w_2) (0) = (1/2, A/4 + \e s_1(1/2, \e), \e s_2(1/2, \e))
\]
which stays on \(M_\e\) and connects the two points \((0, 0, 0)\) and \((1, 0, 0)\). 

Since \( F_\e(w_0) := A(w_0-w_0^2) + \e s_1(w_0, \e) >0\) at \(w_0=1/2\) for sufficiently small \(\e\) and $F_\e(w_0) \to -\infty$ as $\abs{w_0} \to \infty$, there exist $a,b \in \RR$ with \(a < 1/2 < b\) such that \(F_\e(a) = F_\e(b) = 0\) and \(F_\e >0 \) on \((a, b)\). 
Thus, by the succeeding lemma, we have a unique solution \(\wbar_0\colon \RR\to \RR\) of the following ODE 
\begin{equation}\label{one}
w'_0(z) = F_\e(w_0(z)) = A(w_0(z)-w_0(z)^2) + \e s_1(w_0(z), \e)
\end{equation}
subject to the initial datum
\begin{equation}\label{ICone}
w_0(0) = 1/2,
\end{equation}
and the solution is increasing and converges to \(a\) and \(b\) as \(z \to \pm\infty\).
\begin{lemma} \label{lem-odethm}
Let \(B\in \Ccal^1([a, b])\) with \(B(a) = B(b) = 0\),  \(B>0\) on \((a, b)\), and fix \(f_0\in (a, b)\). 
Then, there is a unique \(\Ccal^1\) function \(f\colon \RR\to[a, b]\) satisfying \(f'(z) = B(f(z))\) and \(f(0) = f_0\), which is increasing and \(f(z)\) converges to \(a\) and \(b\) as \(z\to \pm\infty\).    
\end{lemma}
This lemma is a classical result in the theory of ordinary differential equations, specifically in the context of the autonomous case.
For the unique profile $\wbar_0$, we let 
\begin{equation}\label{bw1}
\wbar_1 \coloneqq A(\wbar_0-\wbar_0^2) + \e s_1(\wbar_0, \e),\quad \wbar_2 \coloneqq \e s_2(\wbar_0, \e) \quad\mbox{on } \RR,
\end{equation}
and we will prove that the function \(\wbar \coloneq (\wbar_0, \wbar_1, \wbar_2)\) satisfies \eqref{sys}.

For any fixed \(z_0\in\RR\), the point \(\wbar(z_0)\) is on the manifold \(M_\e\). 
From the Cauchy-Lipschitz theorem, we have a unique Lipschitz solution $\wtil (z) = (\wtil_0, \wtil_1, \wtil_2) (z) $ (locally defined near $z_0$) to the system \eqref{sys} passing through $\wbar(z_0)$ at $z_0$, i.e., $\wtil(z_0)=\wbar(z_0)$.\\ 
Then, since \(\wtil(z_0)=\wbar(z_0) \in M_\e\) by \eqref{bw1} and the local solution $\wtil(z)$ is continuous near $z_0$, it holds from Proposition \ref{prop_FFT} that the local curve $\wtil(z)$ is contained in the locally invariant manifold \(M_\e\), that is,
\beq\label{tils2}
\wtil_1 = A(\wtil_0-\wtil_0^2) + \e s_1(\wtil_0, \e), \quad \wtil_2 = \e s_2(\wtil_0, \e),\quad\mbox{near } z_0.
\eeq
From the definition of \(\wtil\), we have \(\wtil_0' = \wtil_1\), and so \(\wtil_0\) is a local solution of \eqref{one}. 
Then, by the uniqueness of \eqref{one} with $\wtil(z_0)=\wbar(z_0)$, we can deduce that \(\wtil_0 = \wbar_0\) near \(z_0\). 
Furthermore, since both \((\wbar_1,\wbar_2)\) and \((\wtil_1,\wtil_2)\) are parametrized by \(\wbar_0\) and \(\wtil_0\) as in \eqref{bw1} and \eqref{tils2} respectively, those are locally the same as well.
This proves that \(\wbar\) satisfies \eqref{sys} at \(z_0\). 
Since \(z_0\) was chosen arbitrary in \(\RR\), the globally defined function \(\wbar\) is a solution to the system \eqref{sys}. 

Lastly, we prove that \(\wbar\) connects the critical points.
Since \(\wbar_0\) converges at \(\pm\infty\), \(\wbar_1\) and \(\wbar_2\) also converge at \(\pm\infty\) by its constructions.
Hence, the solution \(\wbar\) has a limit at \(\pm\infty\). 
Notice that the limits of a solution of an ODE system at \(\pm\infty\) should be its stationary points.
Thus, \(\wbar\) converges to \((0, 0, 0)\) and \((1, 0, 0)\) at \(\pm\infty\) because there are only two critical points of \eqref{sys} as we observed in advance.

\step{2} In Step 1, for sufficiently small \(\e\), we proved the global existence of \eqref{sys} with
\[
    \lim_{z \to -\infty} w_0(z) = 0, \quad \lim_{z \to \infty} w_0(z) = 1, \quad
    \lim_{z \to \pm\infty} w_0'(z) = 0, \quad \text{ and } \quad \lim_{z \to \pm\infty} w_0''(z) = 0.
\]
Especially, it holds from \eqref{sys} that
\begin{equation}\label{w0ODE}
    \e^2 w'''_0(z) =f(w_0(z), w'_0(z), \e w''_0(z), \e).
\end{equation}
In this step, we will show that the function \(w_0(z)\) in \eqref{w0ODE} gives a solution to the following system:
\begin{equation} \label{intODE}
    \left\{
    \begin{aligned}
        & -\s_\e(\vtil-v_-) - (\util-u_-) = \t \frac{\vtil'}{\vtil}, \\
        & -\s_\e(\util-u_-) + (\ptil-p_-) = \m \frac{\util'}{\vtil}, \\
        & -\s_\e\left(\frac{R}{\g-1}(\thtil-\th_-) + \frac{1}{2}(\util^2-u_-^2)\right)
        + (\ptil\util-p_-u_-) = \k \frac{\thtil'}{\vtil} + \m \frac{\util\util'}{\vtil},
    \end{aligned}
    \right.
    \end{equation}
and thus, to the desired system \eqref{shock_0}.

From now on, we consider the parameter \(\e\) as the amplitude of the shock \(\e = v_+ - v_- > 0\).
First of all, we claim
 \[
    \sbare = \s_\e.
\] 
Given \eqref{end-con}, we have \((p_+-p_-)(v_+-v_-) = -(u_+-u_-)^2\), so that 
\begin{equation}\label{param-th1}
\th_+ = \frac{v_+}{R}\left(p_- - \frac{(u_+-u_-)^2}{v_+-v_-}\right)
= \th_- + \frac{p_-}{R}(v_+-v_-) - \frac{v_+}{R}\frac{(u_+-u_-)^2}{v_+-v_-}.
\end{equation}
Multiplying \eqref{end-con}\({}_{2}\) by \(u_+\), and substracting it from  \eqref{end-con}\({}_{3}\) gives
\[
-\s_\e\frac{R}{\g-1}(\th_+-\th_-) + \frac{\s_\e}{2}(u_+-u_-)^2 + p_-(u_+-u_-) = 0,
\]
which implies 
\begin{equation}\label{param-th2}
\th_+ = \th_- +\frac{\g-1}{2R}(u_+-u_-)^2 - \frac{\g-1}{R}p_-(v_+-v_-).
\end{equation}
From \eqref{param-th1} and \eqref{param-th2}, we obtain the expression for \((u_+-u_-)^2\): 
\begin{equation}\label{param-u}
(u_+-u_-)^2 = \frac{\g p_-(v_+-v_-)}{\frac{\g-1}{2} + \frac{v_+}{v_+-v_-}}
= \frac{\g p_-(v_+-v_-)^2}{v_- + \frac{\g+1}{2}(v_+-v_-)}.
\end{equation}
From \eqref{end-con}\({}_1\), \(u_+ - u_-\) and \(v_+ - v_-\) have the opposite sign. \\
Therefore, we obtain that 
\begin{align} \label{u_+(v_+)}
    u_+ - u_- &= - \frac{\sqrt{\g p_-}}{\sqrt{v_- + \frac{\g+1}{2}(v_+-v_-)}}(v_+-v_-).
\end{align}
Plugging \eqref{u_+(v_+)} into \eqref{end-con}\(_1\), we conclude that 
\begin{equation}\label{s(v_+)s}
\s_\e = -\frac{u_+-u_-}{v_+-v_-} 
= \frac{\sqrt{\g p_-}}{\sqrt{v_- + \frac{\g+1}{2}(v_+-v_-)}}
= \frac{\sqrt{\g p_- }}{\sqrt{v_- + \frac{\g+1}{2}\e}} = \sbare
\end{equation}
as  desired. \\

Now, we define \(\vtil = \vtil(\x)\) by
\begin{equation}\label{vtil}
     \vtil(\x)\coloneqq\e w_0(\e \x)+v_-,\quad \e \x =z.
\end{equation}
Multiplying \(\frac{\e^2 \m\t\k}{R(\e w_0 + v_-)^2}\) on both sides of \eqref{w0ODE} and using \eqref{vtil}, we have
\begin{equation} \label{veqn} 
\begin{aligned}
    0& = \frac{\s_\e\vtil}{\g-1}\left(-\s_\e^2(\vtil-v_-)-\s_\e(\m + \t)\frac{\vtil'}{\vtil}-\m\t\frac{\vtil''}{\vtil^2}+\m\t\frac{(\vtil')^2}{\vtil^3}\right) + \frac{\s_\e p_-}{\g-1}(\vtil-v_-) \\
    &\qquad - \frac{\s_\e^3}{2}(\vtil-v_-)^2 - \s_\e^2\t(\vtil-v_-)\frac{\vtil'}{\vtil} - \frac{\s_\e\t^2}{2}\frac{(\vtil')^2}{\vtil^2} + 
    \s_\e p_-(\vtil-v_-) + p_-\t\frac{\vtil'}{\vtil} \\
    &\qquad + \frac{\k}{R\vtil}\left((p_--2\s_\e^2\vtil+ \s_\e^2v_-)\vtil'-\s_\e(\m+\t)\vtil''-\m\t\left(\frac{\vtil'''}{\vtil}-\frac{3\vtil'\vtil''}{\vtil^2}+\frac{2(\vtil')^3}{\vtil^3}\right)\right)
\end{aligned}
\end{equation}
which can be reduced to
\begin{equation}\label{vODE}
    \begin{aligned}
        &-\frac{\s_\e}{\g-1}\left[ \vtil \left(p_- -\s_\e^2(\vtil-v_-)-\s_\e(\m + \t)\frac{\vtil'}{\vtil}-\m\t\frac{\vtil''}{\vtil^2}+\m\t\frac{(\vtil')^2}{\vtil^3}\right) - p_- v_- \right] \\
        &+\frac{1}{2}\s_\e \left(-\s_\e(\vtil-v_-) - \t\frac{\vtil'}{\vtil} \right)^2 + 
        p_- \left(-\s_\e(\vtil-v_-) - \t\frac{\vtil'}{\vtil}\right) \\
        &\quad = \frac{\k}{R\vtil}\left((p_--2\s_\e^2\vtil+ \s_\e^2v_-)\vtil'-\s_\e(\m+\t)\vtil''-\m\t\left(\frac{\vtil'''}{\vtil}-\frac{3\vtil'\vtil''}{\vtil^2}+\frac{2(\vtil')^3}{\vtil^3}\right)\right).
    \end{aligned}
    \end{equation}
Then we define \(\util\) and \(\thtil\) in accordance with \eqref{intODE}\(_1\) and \eqref{intODE}\(_2\) respectively as follows:
\begin{align} \label{util}
\util\coloneqq u_--\s_\e (\vtil-v_-)-\t \frac{\vtil'}{\vtil}, 
\end{align}
and 
\[
\thtil\coloneqq\frac{p_-}{R}\vtil+\frac{\s_\e}{R}\vtil(\util-u_-)+\frac{\m}{R}\util'.
\]
Since we have
\begin{equation}\label{util'}
\util'= -\s_\e \vtil' - \t \frac{\vtil''}{\vtil} + \t \frac{(\vtil')^2}{\vtil^2},
\end{equation}
\(\thtil\) can be written in the following way:
\begin{equation}\label{th(v)}
    \thtil = \frac{\vtil}{R}\left(p_--\s_\e^2(\vtil-v_-)-\s_\e(\m + \t) \frac{\vtil'}{\vtil} - \m \t \frac{\vtil''}{\vtil^2}+\m \t \frac{(\vtil')^2}{\vtil^3}\right).
    \end{equation}
Hence, differentiating both sides of \eqref{th(v)}, we get
\begin{equation}\label{th'(v)}
    \begin{aligned}
    \thtil' &= \frac{1}{R}\left(p_-\vtil'-2\s_\e^2\vtil\vtil' + \s_\e^2v_-\vtil'-\s_\e (\m + \t)\vtil''- \m \t \left(\frac{\vtil'''}{\vtil}-\frac{\vtil'\vtil''}{\vtil^2} - \frac{2\vtil'\vtil''}{\vtil^2}+\frac{2(\vtil')^3}{\vtil^3}\right)\right) \\
    &= \frac{1}{R}\left((p_--2\s_\e^2\vtil+\s_\e^2v_-)\vtil'-\s_\e(\m+\t)\vtil''- \m \t \left(\frac{\vtil'''}{\vtil}-\frac{3\vtil'\vtil''}{\vtil^2}+\frac{2(\vtil')^3}{\vtil^3}\right)\right).
    \end{aligned}
\end{equation}
Plugging \eqref{util}, \eqref{th(v)}, and \eqref{th'(v)} into \eqref{vODE}, we obtain that
\begin{equation} \label{intODE_1}
    -\s_\e\left(\frac{R}{\g-1}(\thtil-\th_-) - \frac{1}{2}(\util-u_-)^2\right)
    + p_-(\util-u_-) = \k \frac{\thtil'}{\vtil}.
\end{equation}
Then, multiplying \eqref{intODE}\(_2\) by \(\util\), and adding it to \eqref{intODE_1} yields \eqref{intODE}\(_3\):
\[
    -\s_\e\left(\frac{R}{\g-1}(\thtil-\th_-) + \frac{1}{2}(\util^2-u_-^2)\right)
    + (\ptil\util-p_-u_-) = \k \frac{\thtil'}{\vtil} + \m \frac{\util\util'}{\vtil}.
\]
Furthermore, since \(w_0\) is increasing, the monotonicity of \(\vtil\) follows from \eqref{vtil}. 
The monotonicity of \(\util\) and \(\thtil\) follows from \eqref{ratio-vu} and \eqref{ratio-vth}, which will be proved in Section \ref{sec:property}.
\qed

\subsection{Proof of Theorem \ref{thm_EU}: Uniqueness}
In Section \ref{section E}, we proved the existence of a solution of \eqref{intODE} connecting \((v_-, u_-, \th_-)\) to \((v_+, u_+, \th_+)\). Now, we will prove the uniqueness of the monotone solution to \eqref{intODE} connecting two endpoints \((v_-, u_-, \th_-)\) to \((v_+, u_+, \th_+)\).


To achieve this, we will use the following proposition on the unstable manifold at a critical point.

\begin{proposition} \cite[Theorems 9.4, 9.5]{teschl2012ordinary}  \label{TeschlUnique}
Consider a system of ODEs \(z' = F(z)\), where  \(F \in \Ccal^k(\RR^n ; \RR^n)\) and $x_0$ is a critical point, i.e., \(F(x_0) = 0\). Let \(W^-(x_0)\) be the unstable manifold composed of all points converging to \(x_0\) for \(t\to -\infty\). 
Let \(E^+\) (resp. \(E^-\)) denote the linear subspace spanned by eigenvectors of \(D_{x_0} F\) corresponding to negative (resp. positive) eigenvalues. 
Suppose \(D_{x_0} F \) has no eigenvalues on the imaginary axis. \\
Then, there are a neighborhood \(U(x_0) = x_0 + U\) of \(x_0\) and a function \(h^-\in \Ccal^k(E^-\cap U; E^+)\) such that 
\[
W^-(x_0)\cap U(x_0) = \{x_0+a+h^-(a)\mid a\in E^-\cap U\}.
\]
Here, both \(h^-\) and its Jacobian matrix vanishes at \(0\).
\end{proposition}
This implies that the set \(W^-(x_0)\) is a \(C^k\)-manifold of which dimension is equal to the dimension of \(E^-\), and this is tangent to the affine space \(x_0 + E^-\) at \(x_0\).

To apply Proposition \ref{TeschlUnique}, we rewrite the system \eqref{intODE} by using \eqref{intODE_1} as follows:
\begin{equation} \label{intODE2}
    \left\{
    \begin{aligned}
        & \vtil' = \frac{\vtil}{\t}\left[-\s_\e(\vtil-v_-) - (\util-u_-) \right], \\
        & \util' = \frac{\vtil}{\m}\left[-\s_\e(\util-u_-) + (\ptil-p_-)\right], \\
        & \thtil' = \frac{\vtil}{\k}\left[ -\s_\e\left(\frac{R}{\g-1}(\thtil-\th_-) + \frac{1}{2}(\util-u_-)^2 \right) + p_-(\util - u_-)\right].
    \end{aligned}
    \right.
    \end{equation}
It is easy to see from \eqref{end-con} that \((v_\pm, u_\pm, \th_\pm)\) are critical points of the system \eqref{intODE2}. Then we need to compute the sign of eigenvalues of the linearized system at the critical point \((v_-, u_-, \th_-)\)  so as to determine the dimension of the unstable manifold. 
For \eqref{intODE2}, the Jacobian matrix \(J(\vtil, \util, \thtil)\) at \((v_-, u_-, \th_-)\) is given by
\begin{equation}\label{Jacobian}
\left. J(\vtil, \util, \thtil) \right|_{(\vtil, \util, \thtil) = (v_-, u_-, \th_-)}=\begin{pmatrix} 
-\frac{\s_\e v_-}{\t} & -\frac{ v_-}{\t} & 0 \\
-\frac{p_-}{\m} & -\frac{\s_\e v_-}{\m} & \frac{R}{\m} \\ 
0 & \frac{v_-p_-}{\k} & -\frac{R}{\g -1}\frac{\s_\e v_-}{\k}
\end{pmatrix}.
\end{equation}
We then compute the trace and determinant of \eqref{Jacobian} as follows:
\[
\begin{aligned}
\det J((v_-, u_-, \th_-))&= \frac{1}{\m\t\k}\frac{R\s_\e v_-^2}{\g-1}(\g p_--\s_\e^2v_-) > 0, \\
\tr J((v_-, u_-, \th_-))&= -\s_\e v_-\left(\frac{R}{\g-1}\frac{1}{\k}+\frac{1}{\m}+\frac{1}{\t}\right) < 0.
\end{aligned}
\]
Here, we utilized \eqref{s(v_+)2} for the first inequality.
Upon examining the graph of the characteristic equation \(\l^3 + A\l^2 + B\l + C = 0\) with \(C<0\), we deduce that there is at least one positive real root. 
Since the sum of three eigenvalues are negative, the sum of the other two roots is negative and their product is positive; those can be complex roots with negative real part, or two negative real roots.
Hence, the dimension of the unstable manifold, which is equal to the number of eigenvalues with positive real part, is \(1\).

Now, we prove the desired uniqueness. We already constructed a monotone solution \((\vtil, \util, \thtil)\) which satisfies \eqref{intODE} and connects \((v_-, u_-, \th_-)\) to \((v_+, u_+, \th_+)\). 
This solution trajectory \((\vtil, \util, \thtil)\) is locally contained in the one-dimensional unstable manifold \(W^-((v_-, u_-, \th_-))\). More precisely,  
 \(W^-((v_-, u_-, \th_-))\cap \{v\ge v_-, u \le u_-, \th\le \th_-\}\) must coincide with the trajectory \((\vtil, \util, \thtil)\) near \((v_-, u_-, \th_-)\). 
Note that, any other solution which connects \((v_-, u_-, \th_-)\) to \((v_+, u_+, \th_+)\) should be on \(W^-((v_-, u_-, \th_-))\), and especially if it is monotone, it should intersect with the trajectory \((\vtil, \util, \thtil)\).
Hence, by the uniqueness for the Lipschitz autonomous system of ODE, \((\vtil, \util, \thtil)\) is the unique solution in the class of monotone solutions, up to a translation.

\qed

\section{Proof of Theorem \ref{thm_estimates}}
\label{sec:property}
\setcounter{equation}{0}

We present the proof of Theorem \ref{thm_estimates}.
In this section, \(C\) represents a positive constant that may vary from line to line, yet remains independent of \(\e\).

To prove the desired estimates as  in Theorem \ref{thm_estimates}, we will proceed in the following order. First, we obtain \((\abs{\vt''(\x)},\abs{\ut''(\x)},\abs{\tht''(\x)})\le C\e \abs{\vt'(\x)}\). We then prove \eqref{ratio-vu} and \eqref{ratio-vth}, so that \eqref{derivsign} and \eqref{2bound1} are automatically established. Subsequently, we show \eqref{tail-0}, and this with Gronwall's inequality gives \eqref{decay} and \eqref{derivdecay}.

Recall \eqref{sys} and \eqref{vtil}: we have the following relations.
\begin{align} \label{vw12}
\vt'(\x) = \e^2 w'_0(z),
\quad \vt''(\x)  = \e^2 w_2(z) = \e^3 w_1'(z).
\end{align}
Similarly, we also have 
\begin{align} \label{vw34}
\vt'''(\x) = \e^4w''_1(z), \quad \vt''''(\x) = \e^5w'''_1(z).
\end{align}
On the other hand, from \eqref{Me}, we have
\begin{align} \label{w1=}
w_1 = A(w_0-w_0^2) + \e s_1(w_0, \e), \quad w_0\in K,
\end{align}
for some compact set $K$, and  \(s_1\) is smooth.\\
 So, we get \(\abs{w_1} = \abs{w'_0} \le C\). Then, using \eqref{vw12}\(_{1}\), we obtain \(\abs{\vt'(\x)}\le C\e^2\).
Following this, we differentiate \eqref{w1=}:
\[
w'_1 = A(1-2w_0)w'_0 + \e D_1 s_1(w_0, \e)w'_0.
\]
Then, since \(D_1 s_1\) is continuous and \(w_0\) is confined to \(K\), it holds that \(\abs{w'_1(z)}\le C \abs{w'_0(z)}\). 
This together with \eqref{vw12} implies that \(\abs{\vt''(\x)}\le C\e\abs{\vt'(\x)}\). 

To obtain the bounds for \(\abs{\ut''}\) and \(\abs{\tht''}\), we differentiate \eqref{w1=} once and twice more:
\begin{align*}
w''_1 &= -2A (w'_0)^2 + A(1-2w_0)w''_0 + \e D_{11}s_1(w_0, \e)(w'_0)^2 + \e D_1 s_1(w_0, \e)w''_0 \\
&= -2A w_1w'_0 + A(1-2w_0)w'_1 + \e D_{11}s_1(w_0, \e)w_1w'_0 + \e D_1 s_1(w_0, \e)w'_1, \\
w'''_1 &= -6Aw'_0w''_0 + A(1-2w_0)w'''_0 + \e D_{111}s_1(w_0, \e)(w'_0)^3 \\
&\qquad
+ 3\e D_{11}s_1(w_0, \e)w'_0w''_0 + \e D_1 s_1(w_0, \e)w'''_0 \\
&= -6Aw_1w'_1 + A(1-2w_0)w''_1 + \e D_{111}s_1(w_0, \e)w_1^2w'_0 \\
&\qquad
+ 3\e D_{11}s_1(w_0, \e)w_1w'_1 + \e D_1 s_1(w_0, \e)w''_1.
\end{align*}
Then, since all the higher derivatives of \(s_1\) are continuous and \(w_0\) is confined to \(K\) as well, we have \(\abs{w''_1(z)}\le C\abs{w'_0(z)}\) and \(\abs{w'''_1(z)}\le C\abs{w'_0(z)}\).\\
Hence, from \eqref{vw34}, we also obtain \(\abs{\vt'''(\x)}\le C\e^2\abs{\vt'(\x)}\) and \(\abs{\vt''''(\x)}\le C\e^3\abs{\vt'(\x)}\). \\
Thus, by differentiating \eqref{util'} and \eqref{th'(v)}, we use all of these estimates to yield that 
\[
\abs{\ut''(\x)}, \abs{\tht''(\x)} \le C\e\abs{\vt'(\x)}.
\]
We now show \eqref{ratio-vu} and \eqref{ratio-vth}: from \eqref{shock_0}\(_1\), it follows that
\[
\abs{\ut' + \s_\e\vt'} = \abs{\t \frac{\vt''}{\vt}-\t\frac{(\vt')^2}{\vt^2}} \le C\e\abs{\vt'}.
\]
Then, \eqref{sm1} gives the desired result \eqref{ratio-vu}.
Likewise, for \eqref{ratio-vth}, we use \eqref{th'(v)} to get
\[
\abs{\tht' - \frac{1}{R}((p_--2\s_\e^2\vt+ \s_\e^2v_-)\vt')} \le C\e \abs{\vt'}.
\]
Since \(\vt\) is increasing, we have \(\abs{\vt-v_-} \le \abs{v_+-v_-} = \e\). Considering this and \eqref{sm1}, the definition of \(\s_*\) gives the desired one \eqref{ratio-vth}.
Furthermore, these two estimates with \(\vt'>0\) imply that \(\abs{\vt'}\sim \abs{\ut'}\sim |\tht'|\) and \eqref{derivsign}.

In order to prove the exponential decays \eqref{decay} and \eqref{derivdecay}, it is required to show \eqref{tail-0} first. 
To this end, we rewrite \eqref{vODE} in the following form:
\begin{align}
\begin{aligned} \label{vODEre}
&\frac{\s_\e^3 \vt}{\g-1}(\vt-v_-) - \frac{\s_\e p_-}{\g-1}(\vt-v_-) + \frac{\s_\e^3}{2}(\vt-v_-)^2 - \s_\e p_-(\vt-v_-) \\
&= \frac{\s_\e\vt}{\g-1}\left(-\s_\e(\m + \t)\frac{\vt'}{\vt}-\m\t\frac{\vt''}{\vt^2}+\m\t\frac{(\vt')^2}{\vt^3}\right) 
- \s_\e^2\t(\vt-v_-)\frac{\vt'}{\vt} - \frac{\s_\e\t^2}{2}\frac{(\vt')^2}{\vt^2} + p_-\t\frac{\vt'}{\vt} \\
&\qquad + \frac{\k}{R\vt}\left((p_--2\s_\e^2\vt+ \s_\e^2v_-)\vt'-\s_\e(\m+\t)\vt''-\m\t\left(\frac{\vt'''}{\vt}-\frac{3\vt'\vt''}{\vt^2}+\frac{2(\vt')^3}{\vt^3}\right)\right).
\end{aligned}
\end{align}
Afterwards, we analyze both sides of \eqref{vODEre} respectively.
For the left-hand side, we use \eqref{s(v_+)} to have the following:
\begin{equation}
\begin{aligned}\label{tail-LHS}
LHS
&= \s_\e^3(\vt-v_-)\left[\frac{-\g p_-}{\s_\e^2 (\g-1)} + \left(\frac{\vt}{\g-1} + \frac{1}{2}(\vt-v_-)\right)\right] \\
&= \s_\e^3(\vt-v_-)\left(-\frac{v_+ + \frac{\g-1}{2}(v_+-v_-)}{\g-1}+\frac{\vt}{\g-1} + \frac{1}{2}(\vt-v_-)\right) \\
&= \s_\e^3\left(\frac{1}{\g-1} + \frac{1}{2}\right)(\vt-v_-)(\vt-v_+).
\end{aligned}
\end{equation}
We then observe that all the terms on the right-hand side can be controlled by \(\vt'\). Indeed,
\begin{equation}
\begin{aligned}\label{tail-RHS}
RHS
&= \left(-\frac{\s_\e^2(\m+\t)}{\g-1} + p_-\t\frac{1}{\vt} -2\frac{\k}{R}\s_\e^2 + \frac{\k}{R\vt}(p_- + \s_\e^2 v_-) + \Ocal(\e)\right)\vt' \\
&= \left(-\frac{\g(\m+\t)p_-}{(\g-1)v_-} + p_-\t\frac{1}{v_-} -\frac{(\g-1)\k p_-}{Rv_-} + \Ocal(\e)\right)\vt' \\
&= -\frac{p_-}{v_-}\left(\frac{\g}{\g-1}\m + \frac{1}{\g-1}\t + \frac{\g-1}{R}\k + \Ocal(\e)\right)\vt'.
\end{aligned}
\end{equation}
Combining \eqref{vODEre} with \eqref{tail-LHS} and \eqref{tail-RHS}, we obtain \eqref{tail-0}, i.e., 
\[
C^{-1}(v_+-\vt)(\vt-v_-) \le \vt' \le C(v_+-\vt)(\vt-v_-).
\]
With this, the standard argument of Gronwall's inequality (for example, see \cite[Lemma2.1]{KV21}) gives that
\begin{equation}
\begin{aligned}\label{vdecay}
\abs{\vt(\x)-v_-}\le C\e e^{-C\e\abs{\x}} \text{ for  } \x \le 0, \quad\quad
\abs{\vt(\x)-v_+}\le C\e e^{-C\e\abs{\x}} \text{ for  } \x \ge 0.
\end{aligned}
\end{equation}
Then, by plugging these estimates into \eqref{tail-0}, we get
\[
\abs{\vt'(\x)} \le C\e^2e^{-C\e\abs{\x}},
\]
and this together with \eqref{ratio-vu} and \eqref{ratio-vth} gives the full strength of \eqref{derivdecay}.\\
In addition, applying \eqref{vdecay} and \eqref{derivdecay} to \eqref{intODE}\({}_1\) and \eqref{intODE_1}, we obtain that 
\begin{align*}
&\abs{\ut(\x)-u_-}\le C\e e^{-C\e\abs{\x}} \text{ for  } \x \le 0, 
&&\abs{\ut(\x)-u_+}\le C\e e^{-C\e\abs{\x}} \text{ for  } \x \ge 0, \\
&\abs{\tht(\x)-\th_-}\le C\e e^{-C\e\abs{\x}} \text{ for  } \x \le 0, 
&&\abs{\tht(\x)-\th_+}\le C\e e^{-C\e\abs{\x}} \text{ for  } \x \ge 0.
\end{align*}
This completes the proof of \eqref{decay} and Theorem \ref{thm_estimates}.
\qed

\bibliographystyle{plain}
\bibliography{reference} 

\end{document}